\documentclass[twoside,12pt, leqno]{amsart}
\usepackage{amsmath,amscd,amsthm,amssymb,amsxtra,latexsym,epsfig,epic,graphics}
\usepackage[matrix,arrow,curve]{xy}
\usepackage{graphicx}
\usepackage{svg}
\usepackage{diagrams,enumitem}
\usepackage{tikz}  
\usepackage{color} 
\usetikzlibrary{arrows,calc} 
\usetikzlibrary{decorations.pathmorphing}
\usetikzlibrary{decorations.markings} 
\usetikzlibrary{decorations.pathreplacing} 
\usetikzlibrary{plothandlers}

\def\antiddot{\mathinner{\mkern1mu\raise1pt\vbox{\kern7pt\hbox{.}}\mkern2mu
        \raise4pt\hbox{.}\mkern2mu\raise7pt\hbox{.}\mkern1mu}}

\newcommand{\FF}{{\mathbb F}}
\newcommand{\GG}{{\mathbb G}}

\newcommand{\lcm}{{\rm{lcm}}}

\newcommand{\Ext}{{\rm{Ext}}}

\newcommand{\punkt}{\hspace{-.3ex}\raise.15ex\hbox to1ex{\Huge.}}

\def \fix#1 {{\hfill\break \bf (( #1 ))\hfill\break}}

\DeclareMathOperator{\reg}{reg}

\DeclareMathOperator{\Hom}{Hom}

\DeclareMathOperator{\im}{im}

\DeclareMathOperator{\Tor}{Tor}

\newtheorem{theorem}{Theorem}[section]
\newtheorem{lemma}[theorem]{Lemma}
\newtheorem{proposition}[theorem]{Proposition}
\newtheorem{corollary}[theorem]{Corollary}

\theoremstyle{definition}
\newtheorem{definition}[theorem]{Definition}
\newtheorem{question}[theorem]{Question}

\newtheorem{remark}[theorem]{Remark}

\newtheorem{example}[theorem]{Example}

\def\FF{{\mathbb F}}

\def\fix#1{{\bf ***Fix:} #1 {\bf ***}}

\def\mm{{\frak m}}

\def\m{{\frak m}}

\def\pol{{\rm pol}}
\def\depol{{\rm depol}}

\def\lbracket{{[\kern-1.5pt[}}
\def\rbracket{{]\kern-1.5pt]}}

\makeatletter
\def\Ddots{\mathinner{\mkern1mu\raise\p@
\vbox{\kern7\p@\hbox{.}}\mkern2mu
\raise4\p@\hbox{.}\mkern2mu\raise7\p@\hbox{.}\mkern1mu}}
\makeatother

\usepackage{times}
\newdimen\x \x=12pt

\usepackage{color}

\usepackage[breaklinks,bookmarksopen,bookmarksnumbered,urlcolor=blue]{hyperref}
\hypersetup{colorlinks=true,backref=true,citecolor=blue}
\usepackage{cleveref}
\author[Hailong Dao]{Hailong Dao}
\address{Department of Mathematics, University of Kansas, Lawrence, KS 66045-7523, USA}
\email{hdao@ku.edu}

\author[David Eisenbud]{David Eisenbud}
\address{Department of Mathematics, University of California at Berkeley and the Mathematical
Sciences Research Institute, Berkeley, CA 94720, USA}
\email{de@msri.org}


\title{Linearity of Free Resolutions of Monomial Ideals}
\begin{document}
\dedicatory{Dedicated to Juergen Herzog, inspiring mathematician and master of monomials,\\on the occasion of his 80th birthday!}

\begin{abstract}
We study monomial ideals with linear presentation or partially linear resolution. We give combinatorial characterizations of linear presentation for square-free ideals of degree 3, and for primary ideals whose resolutions are linear except for the last step (the ``almost linear'' case). We also give sharp bounds on Castelnuovo-Mumford regularity and numbers of generators in some cases. 

It is a basic observation that linearity properties are inherited by the restriction of an ideal to a subset of variables, and we study when the converse holds. We construct fractal examples of almost linear primary ideals with relatively few generators related to the Sierpi\'nski triangle. Our results also lead to classes of  highly connected simplicial complexes $\Delta$ that can not be extended to the complete $\dim \Delta$-skeleton of the simplex on the same variables by shelling.    
\end{abstract}

\maketitle

\section*{Introduction} 
Fix an ambient dimension $n$ and a degree $d$.
Let $S = k[x_{1},\dots,x_{n}]$ be a polynomial ring over a field $k$, and set $\mm = (x_1,\dots,x_n)$. For any finitely generated graded $S$-module  $M$ we write
$$
t_s(M) :=  \max\{e \mid \Tor_s^S(M,k)_e \neq 0\}.
$$
 We will use these definitions throughout the paper.

We say that a homogeneous ideal $I\subset S$ satisfies the condition $N_{d,p}$ if $t_s(I) = d+s$ for all $s\leq p-1$. Thus $N_{d,1}$ is the condition that $I$ is generated in degree $d$, $N_{d,2}$ adds the condition that $I$ is linearly presented, and more generally $N_{d,p}$ is the condition that $I$ has a linear resolution for $p-1$ steps. Green's condition $N_p$ is, in this notation, $N_{2,p}$. We describe an ideal $I$ satisfying $N_{d,q}$ as having \emph{linear resolution}, here $q$ is the projective dimension of $S/I$, and that an ideal satisfying $N_{d,q-1}$ as having \emph{almost linear resolution}.

Studying ideals satisfying $N_{d,p}$ is the same as studying the successive maxima of $t_s(I)$ for arbitrary ideals $I$:

\begin{proposition}\label{trunc}(Truncation principle, \cite[Proposition 1.7]{EHU}) Let $I$ be a homogeneous ideal of $S$ as above. For any integer $s\geq 0$,
the ideal $J = I\cap \mm^{t_s(I)-s}$ has linear resolution for $s$ steps, while for $r \geq s$ we have $t_r(J) = t_r(I)$; thus
$I\cap \mm^d$ satisfies $N_{d,p}$ for $p = \max\{a \mid t_a(I) \leq d+a\}+1$.
 \end{proposition}

In this paper we focus on monomial ideals.
The square-free monomial ideals satisfying $N_{2,n}$ were classified in a famous paper of Fr\"oberg~\cite{Fr}, and the result was extended to a description
of monomial ideals satisfying any $N_{2,p}$ in \cite[Theorem 2.1]{EGHP}. However, concrete characterizations of  monomial ideals satisfying $N_{d,p}$ for $d\geq 3$ are unknown in general, and many basic questions about them have not been thoroughly investigated. Can we characterize them combinatorially? What is the computational cost of checking whether an ideal is $N_{d,p}$? Are there sharp bounds on Betti numbers of these ideals, in particular their number of generators and regularity? What about those that achieve such bounds? See section \ref{sec_ques} for a more detailed discussion of these questions together with brief reviews of relevant literature. 

Work of Boocher and Peeva-Velasco establishes a \emph{locality principle}: the condition $N_{d,p}$ is inherited by the ideals generated by various subsets of generators. The consequences of this, worked out in section \ref{sec_prin}, are used throughout this paper. Conversely, if the restriction of a monomial ideal to sufficiently large subsets of the variables satisfies $N_{d,p}$, then the same is true of the whole ideal. For example, linear presentation can be checked by restricting to 2d variables; but for cubic ideals,
4 variables (plus an auxiliary condition) is enough, as we establish in \ref{sec_linear}. It would be interesting to know optimal results of this type more generally.

Our best results concern primary monomial ideals. We give sharp regularity bounds for such ideals that satisfy $N_{d,p}$ in \ref{sec_reg}. It is well-known that $\mm^d$ is the only $\mm$-primary ideal with linear resolution. We give a constructive characterization of $\mm$-primary ideals with almost linear resolution---that is, satisfying $N_{d,n-1}$.
We also show that a primary cubic monomial ideal with linear presentation must contain the degree 3 part of the ideal generated by the squares of variables; and given that condition, linear presentation can be tested by restricting to monomials in just 4 of the variables at a time. 
A fractal construction in  \ref{sec_last}, related to the Sierpin\'nski triangle, yields such ideals whose number of generators is an arbitrarily small fraction of the number of generators of $\mm^d$. 

An old question (both for square-free and other ideals) asks, given a monomial ideal satisfying some $N_{d,p}$, when can one adjoin one more monomial, keeping the linearity?
Using our structure theory for almost linear primary ideals, we give many examples where this is not possible. For instance, if a primary monomial ideal satisfying $N_{d,n-1}$ has regularity at least $d+2$, then adding a monomial can never both preserve linear presentation and also change the regularity---thus, for example, it is never possible
to reach $\mm^d$, an ideal of regularity $d$, by adding one monomial at a time while preserving linear presentation. Polarizing such examples, we obtain square-free examples as well. These square-free monomial ideals correspond to
examples of highly connected simplicial complexes (i.e, satisfying Serre's condition $(S_l)$) that can't be extended to the full skeleton of the simplex on all variables using shelling moves.  

We collect and discuss some of our favorite open questions in \cref{sec_ques}.\\

\noindent\textbf{Acknowledgements}: The first author is grateful to the
National Science Foundation for partial support. The second author acknowledges support from a Simons Collaborator Grant and from University of Kansas.
The authors would like to thank Kangjin Han, Siamak Yassemi and the referee for helpful comments and corrections that improve the readability of the paper.\\

\noindent\textbf{Data Availability Statement}: Data sharing is not applicable to this article as no new data were created or analyzed in this study.

\section{Locality}\label{sec_prin}

Let $I\subset S$ be a homogeneous ideal.  If $I$ is a monomial ideal,  and $m$ is a monomial, let $I_{\leq m}$ denote the ideal generated by monomial generators of $I$ that divide $m$. If $K$ is a subset of $\{x_1,\dots, x_n\}$ we write $I_K$ for the ideal obtained from $I$ by  restricting to the variables in $K$ (i.e., setting all the variables not in $K$ to zero).  
We begin by applying a result of Peeva and Velasco that extend work of Boocher to describe the minimal free resolution of $I_\leq m$ as a subcomplex of the minimal
free resolution of $I$:

\begin{theorem}\cite[Proposition 3.10]{PV} \label{pv}
If $F$ is the multi-graded minimal free resolution of a monomial ideal $I\subset S$, and $m\in S$ is a monomial, then the minimal free resolution of $I_{\leq m}$ is the subcomplex
of $F$ formed from all summands of terms in $F$ whose degree divides $m$.\qed
\end{theorem}

A first consequence is that we can make many ideals satisfying $N_{d,p}$ from one of them:

\begin{corollary}\label{leqm}
If $I$ is a monomial $N_{d,p}$ ideal, then so is $I_{\leq m}$ for any monomial $m$. In particular, 
\begin{enumerate}[label=(\alph*)]
\item The square-free part of $I$ (the ideal generated by square-free monomial generators of $I$) is also $N_{d,p}$.
\item The restriction of $I$ to any $r\leq n$ variables is also $N_{d,p}$. 
\end{enumerate}
\end{corollary}

\begin{proof}
The first assertion is immediate from~\ref{pv}. For $a)$, take $m$ to be the product of all variables. For $b)$, harmlessly supposing that the variables are $x_1,\dots, x_r$, we take $m=x_1^d\dots x_r^d$. 
\end{proof}

Theorem~\ref{pv} yields a \emph{locality principle}: the condition $N_{d,p}$ is determined by relatively small subsets of the generators of the ideal, and by the restrictions to relatively few variables. In the (generally nonminimal) Taylor resolution $G$ of a monomial ideal $I$, the degrees of the generators of $G_{s}$ are the least common multiples of $s+1$ minimal generators of $I$, and thus these are the only degrees that can occur among generators of the $s$-th module
in a minimal free resolution.  Combining this Theorem~\ref{pv} with the behavior of the Taylor resolution we deduce:


\begin{corollary}\label{local} (Locality principle)
Let $I\subset S$ be a  monomial ideal.  If $F$ is the multi-graded minimal free resolution of $I$, then the degrees of the homogeneous generators of $F_s$  also appear in the minimal free resolution  of an ideal $I_{\leq m}$ for some
monomial $m$ that is the least common multiple of $s+1$ minimal generators of $I$.

Thus the following are equivalent:
\begin{enumerate}
\item $I$ satisfies $N_{d,p}$
 \item  $I_{\leq m}$  satisfies $N_{d,p}$ for all
least common multiples $m$ of $p$ of the monomial generators of $I$.
\item The restriction of $I$ to any $r=dp$ variables satisfies $N_{d,p}$. 
\end{enumerate}
\end{corollary}

%
%

\begin{proof}
(2) follows from (3) since the $\lcm$ of $p$ generators involves at most $dp$ variables.
\end{proof}

This result implies that $\mm$-primary ideals satisfying $N_{d,p}$ cannot be too small. Write $\mm^{[t]}$ for the ideal $(x_1^t,\dots, x_n^t)$. 

\begin{theorem}\label{ndd1}
Let $I$ be a $\mm$-primary monomial ideal satisfying $N_{d,p}$. The following hold: 
\begin{enumerate}
\item $I$ contains   $\sum_{\{i_1,\dots, i_p\}\subset [n]}(x_{i_1},\dots, x_{i_p})^d$.
\item If $p\geq \min\{n,d\}$, then $I=\mm^d$. 
\item If $d\geq p$ then $I$ contains $\mm^{[d-p+1]}\mm^{p-1}$. 
\end{enumerate}
\end{theorem}

\begin{proof}
(1) By the locality principle, it suffices to prove that if the number of variables $n$ is equal to $p$, then $I = \mm^{d}$, and this is well-known:
In this case $I$ has linear resolution, so $S/I$ has regularity $d-1$. Since $I$ is $\mm$-primary monomial ideal, this implies that $S/I$ is zero in degrees $\geq d$---that is, $\mm^{d}\subset I$.

(2) The case $p\geq n$ follows from (1). 
On the other hand, a monomial of degree $d$ can contain at most $d$ variables, so if $p\geq d$, then (1) shows that $I$ contains every monomial of degree $d$.

(3) We have to show that every degree $d$ monomial of the form $m = x_{i}^{d-p+1}m'$ is in $I$. Since $m'$ has degree $p-1$,
it can be divisible by at most $p-1$ variables, so $m$ is divisible by at most $p$ variables. By (1) we have $m\in I$.
\end{proof}

\section{Linearly presented monomial ideals}\label{sec_linear}
\begin{definition}\label{defdual}
Let $I$ be a monomial ideal. We define the dual graph of $I$, $G(I)$ as follows: the vertices of $G(I)$ are the minimal monomial generators of $I$, and there is an edge between $f,g$ if and only if $|\gcd(f,g)| = |f|-1 = |g|-1$ (equivalently $|\lcm(f,g)| = |f|+1 = |g|+1$).   
\end{definition}

The following is well-known to experts (see for instance \cite[Proposition 2.1, Corollary 2.2]{BHZ}), we include it here with a short proof for the convenience of the readers. 
\begin{proposition}\label{lppath}
A monomial ideal $I$ is linearly presented (i.e, is $N_{d,2}$) if and only if $G(I_{\leq m})$ is connected for $m=\lcm(f,g)$ where $f,g$ are any  two minimal monomial generators  of $I$. More concretely, the condition $N_{d,2}$ of $I$ is equivalent to the following: given any monomial generators $f,g$ of $I$, there is a path connecting $f,g$ whose vertices are generators dividing $\lcm(f,g)$. 
\end{proposition}
\begin{proof}
By Corollary \ref{local}, $I$ is linearly presented if and only if any $I_{\leq m}$ is linearly presented for such $m$. By the formula computing Betti numbers for monomial ideals using $\lcm$ lattices (\cite[Theorem 2.1]{GPW}), this is equivalent to the open interval below $m$ being connected for any $m$ of size at least $d+2$, which is equivalent to $G(I_{\leq m})$ being connected for all $m$. 
\end{proof}

The usual characterization of quadratic square-free ideals $I$ satisfying $N_{2,2}$ is that the 1-skeleton of the Stanley-Reisner simplicial complex associated to $I$ should have no induced cycle of length 4 without a chord; and this comes down to saying that the restriction of $I$ to 4 variables cannot be $x_{1}x_{2}, \ x_{3}x_{4}$. By Corollary~\ref{local} (3), the condition $N_{3,2}$ can be decided by the restrictions of $I$ to subsets of 6 variables. We have the following characterization:

\begin{theorem}\label{cubic linear pres}
A square-free monomial ideal  $I$  generated in degree $3$ is linearly presented if and only if the restriction of  $I$  to at most $6$ variables is not, up to relabeling of variables, a disconnected (in the sense of the dual graph described in \ref{defdual} and \ref{lppath}) subset of either
\begin{enumerate}
\item ${x_1x_2x_3,\ x_4x_5x_6}$; or
\item ${x_1x_2x_3,\ x_1x_2x_4,\ x_1x_2x_5,\ x_3x_4x_5}$.
\end{enumerate}
\end{theorem}

A similar result is announced in \cite[Theorem 2.2]{FSY}.

\begin{proof}
If the restriction is disconnected then, by Proposition~\ref{lppath}, $I$ is not linearly presented. 

Conversely, if $I$  is not linearly presented then there is a pair of generators $f, g \in I$ such that there is
no path from $f$ to $g$ within the monomials supported in the support of $fg$. If the support of $fg$ were just $4$ variables this would be impossible. If the support of $fg$ is 5 variables, then we must show we are in case (2). (writing numbers in place of variables for clarity) we can assume that $f,g$ are $123$ and $345$. If there were another generator of $I$
with support in $12345$ it cannot contain 3, since then it would form a path from $f$ to $g$. Thus we may assume
that it is $124$. Now the only additional monomials that could be in $I$ without forming a path would be those
in (2).

Finally, 
If $f,g$ involve $6$ variables we may assume that they are $123$ and $456$. If there is no other generator then we are in case (1). Otherwise, there is another generator with those variables, and we may suppose that it is $345$, which is directly
connected to $456$. Thus the restriction to $12345$ must be disconnected, and we are in case (2).
\end{proof}

In the case of a primary ideal generated by cubics, we can do with restrictions to fewer than 6 variables:

\begin{theorem}\label{primecubic}
Let $I$ be a $\mm$-primary monomial ideal generated in degree $3$. Then $I$ is linearly presented if and only if the following hold:
\begin{enumerate}
\item $I$ contains $\mm^{[2]}\mm$ (in other words $I$ contains all non-square-free cubics). 
\item The restriction of $I$ to any four distinct variables contains  at least  two square-free cubics or none.  
\end{enumerate}

Consequently, a primary ideal generated by cubic monomials is linearly presented if and only if its restriction to any four variables is linearly presented. 
\end{theorem}

\begin{proof}
The necessity of  $(1)$ follows from Theorem \ref{ndd1}. If $I$ is linearly presented, then so is $I$ restricted to four variables, say $J := I_{\{a,b,c,d\}}$.  If $J$ contains only one  square-free cubic, say $abc$, then there is no path between $abc$ and $ad^2$, so $J$ is not linearly presented. This shows the necessity of (2). 

Conversely, suppose that $I$ satisfies $(1)$ and $(2)$. To prove that $I$ is $N_{3,2}$, it is enough to check the connectivity condition of Proposition \ref{lppath}. Let $I'=\mm^{[2]}\mm$.  Let $f, g \in I$ be generators. If they are both in $I'$, then since $I'$ is linearly presented (use \ref{trunc}), we know that there is a path between them in $G(I')$, and hence also in $G(I)$. So we can assume one of them is square-free, say $f=abc$. There are now several cases up to permutations. 
If the degree of the $\lcm$ of $f,g$ is $4$, they are directly connected, so we may assume that the degree $\ell$
of the  $\lcm$ is 5 or 6. Up to permutation of variables we may assume that:
\begin{itemize}
 \item If  $\ell = 5$, then $g$ is one of
 $$
c^3, c^2d,  cd^2, cde, 
$$
\item If $\ell =6$, then $g$ is one of
$$
d^3, d^2e, def, 
$$
\end{itemize}
and in each case we must construct a sequence of monomials in $I$ starting with $f$ and ending with $g$ such that consecutive pairs have lcm of degree 4 and all the elements
divide $\lcm(f,g)$. We give a suitable path for each case:

\noindent $g = c^3$:\quad $abc, bc^2, c^3$ satisfies the hypothesis because $b^2c\in I$ by condition (1).\\
\noindent $g = c^2d$: \quad $abc, bc^2,c^2d$.\\
\noindent $g = cd^2$: \quad By condition (2) there must be another square-free monomial in $I$ that divides $abcd$, and it must be divisible by $d$; up to permutation
it is say $abd$ or $acd$. In the first case we have the path $abc, abd, ad^2, cd^2$, while in the second we have $abc, acd,cd^2$.\\
\noindent $g = cde$: \quad Again either $abd$ or $acd$ is in $I$. Restricting to $(a,b,c,d)$, $(2)$ tells us there is a square-free cubic other than $abc$. If that cubic is directly connected to $cde$ we are done; in the contrary case it must be $abd$. Similarly starting from $cde$ we may assume, $acd \in I$, so we can take the path $abc,abd,acd, cde$.\\
\noindent $g = d^3$: \quad By (2) we may assume that $abd\in I$, so we have the path $abc,abd,ad^2, d^3$.\\
\noindent $g = d^2e$: \quad Starting as for $d^3$ we get the path $abc,abd,ad^2, d^2e$.\\
\noindent $g = def$: \quad Any of the possible paths from $abc$ to $cde$ considered in the case $g = cde$ extends to $def$.

The last assertion follows because (1) and (2) can be checked by restricting to at most four variables.

\end{proof}

%
%
%
%
%

\begin{example}
The size ($4$) of subset of variables needed to test linear presentation in \ref{primecubic} is optimal. Consider the ideal $I=(a^2,b^2,c^2,d^2)(a,b,c,d)+(abc)$. $I$ is primary, and its restriction to any three variables is linearly presented, but $I$ itself is not. 
\end{example}

\section{Regularity bounds for $N_{d,p}$ ideals}\label{sec_reg}

If $M$ is a finitely-generated graded $S$-module the (Castelnuovo-Mumford) \emph{regularity} of $M$ is defined to be
 $\reg M := \max_s \{t_s(M)-s\}.$ There has been considerable interest in bounding the regularity
under various assumptions on $M$. It turns out that the bound for $\mm$-primary monomial ideals is much smaller than that for non-monomial ideals, which was given
in \cite{MN}.

\begin{theorem}\label {regularity bound}
Suppose that $I\subset S$ is an $\mm$-primary ideal satisfying $N_{d,p}$ with $p\geq 1$.
\begin{enumerate}
 \item If $p =n - 1$ then $\reg I \leq 2d-1$, and this bound is sharp for all $n,d$.
 \item If $I$ is generated by monomials and $p\leq n$ then 
 $$
 \reg(I) \leq d + (n-p) \lfloor \frac{d-1}{p}\rfloor\,.
 $$
 This bound is sharp for all $n,d,p$ with $p\leq \min\{n,d\}$.
 In particular,
if $p=n-1$ then $\reg I \leq d+\lfloor \frac{d-1}{n-1}\rfloor\, .$
 \end{enumerate}
\end{theorem}

Note that if $I$ is a monomial ideal satisying $N_{d,p}$ with $p\geq \min\{n,d\}$ then $I = \mm^d$ by Theorem~\ref{ndd1}(2), so $\reg I = d$.

\begin{proof}[Proof of Theorem \ref{regularity bound}]
The inequality in item (1) follows from formula (1) in Section 11 of \cite{EHU}. 

For the inequality in item (2), let $m_1,\dots, m_k$ be minimal monomial generators of $I:\mm$. The regularity of $I$
is one more than the maximum of the degrees of the $m_i$, and $I$ is generated by the monomials
of degree $d$ that do not divide any of the $m_i$.  

By Theorem~\ref{ndd1}, $I$ contains every monomial of degree $d$ that involves only $p$ variables, so no $m_i$  is divisible by a monomial of degree $d$ in just $p$ variables; that is,
the sum of the largest $p$ exponents of $m_i$ is at most $d-1$. If we order the variables so that the exponent of $x_j$ in $m_i$ is a non-increasing function of $j$, then the maximum possible degree of $m_i$ is achieved
if the sum of the first $p$ exponents is $d-1$, and the rest of the exponents are equal to the $p$-th exponent.
The largest value that the $p$-th exponent could have is $\lfloor(d-1)/p\rfloor$. Thus 
$$
1+\deg m_i \leq d + (n-p) \lfloor \frac{d-1}{p}\rfloor\,
$$
proving the inequality. 

To complete the proof, we give examples of
ideals that achieve the bounds.

\def\F{{\mathbb F}}
\begin{example}\label{sharp examples}
(1) To see that the bound in (1) is sharp, suppose that the field $k$ has characteristic 0, and set $A = S/(x_1^d,\dots, x_n^d)$, so that the socle of $A$
is generated by $\prod_{i=1}^nx_i^{d-1}$, which has degree $n(d-1)$. The element $\sigma = \sum_i x_i$ is a strong Lefschetz element for $A$ (see, for instance, \cite[Theorem 1.1]{MN}); that is, multiplication by a power of $\sigma$ induces an isomorphism $A_e \to A_{n(d-1)-e}$ for every $e$. Set
$$
I=(x_1^d,\dots, x_n^d):\sigma^{(n-2)(d-1)} 
$$
Since $A$ is Gorenstein, so is $S/I$, and the socle of $S/I$ is in degree $n(d-1) - (n-2)(d-1) = 2d-2$, so $\reg I = 2d-2$. Since multiplication by $\sigma^{(n-2)(d-1)}$ induces an isomorphism from $A_{d-1}$
to $A_{(n-1)(d-1)}$, we see that $I\subset \mm^d$. Let $\F$ be the minimal $S$-free resolution of $I$. It follows that the $j$-th term
$F_j$ of $\FF$ is generated in degrees $\geq d+j-1$, and since the $S$-free resolution of $S/I$ is symmetric this must be an equality for $j<n$. Thus $I$ satisfies $N_{d,n-1}$.

(2) To see that the bound in (2) is sharp, let $q = \lfloor (d-1)/p\rfloor$ and write $d-1 = qp + r$, with $r < p$.  Set $m = (x_1\cdots x_n)^qu$, where $u$ is any monomial of degree $r$, so that $\deg m = qn+r$. Let $a_i$ be the exponent of $x_i$ in $m$. Reorder the variables if necessary
so that $a_1\geq\cdots \geq a_n$; note that $a_p = \cdots = a_n$ since $r<p$. 

Set $J = (x_1^{a_1+1}, \dots, x_n^{a_n+1})$ and note that $J:\mm= J+(m)$.
Finally, let
$I = \mm^d\cap J$. Since $a_i\leq d$ for all $i$, the ideal $I$ is generated by the monomials of degree $d$ not dividing $m$, and $I:\mm = I+(m)$.

By Proposition~\ref{trunc} the ideal
$I$ satisfies $N_{d,p}$. We have
\begin{align*}
 \reg I &= 1+\deg m\\
 &= 1+a_1+\cdots +a_n\\ 
 &= 1+a_1+\cdots+a_p + (n-p)q = d-1+(n-p) \lfloor \frac{d-1}{p}\rfloor\,,
\end{align*}

as required.
%
%
%
%
\end{example}
\end{proof}

\section{Almost linear resolutions: the condition $N_{d,n-1}$}\label{sec_almost}

In this section we give a characterization---in some sense a parametrization---of primary monomial ideals satisfying $N_{d,n-1}$, that is, with almost linear resolution. We will state the condition in terms of three definitions:

\begin{definition} 
We define the \emph{$s$-shadow} of a set of monomials $m_1,\dots, m_u$ to be the set of all monomials of degree $s$ that divide some $m_i$.

We say that monomials $m,m'$ are \emph{$s$-separated} if $\deg \gcd(m,m') < s$, or equivalently if their $s$-shadows do not intersect.

We say that a monomial $m$ is \emph{$s$-saturated} if $m$ is divisible by every monomial of degree $\leq \deg m - s$, or equivalently, if 
the exponent of each variable in $m$ is at least $\deg m - s$.

\end{definition}

\begin{theorem}\label{n-1}
 Suppose that $I\subset k[x_{1},\dots, x_{n}]$ is an $\mm$-primary monomial ideal generated in degree $d$, and let $N$ be the set of monomials of degree $d$ that
 are \emph{not} in $I$. 
 
 The ideal $I$ satisfies $N_{d,n-1}$ if and only if $N$ is the shadow of a set $\{m_1,\dots, m_u\}$ of
  $(d-1)$-saturated monomials
  that are pairwise $(d-1)$-separated. In this case $m_1,\dots, m_u$ are generators of the socle of $S/I$, and thus $\reg I = 1+\max_i \deg m_i$. 
\end{theorem}

\begin{example}\label{ex3}
When $n = 3$ the condition $N_{d,n-1} = N_{d,2}$ is the condition that the ideal $I$ is linearly presented. The set of monomials of a given degree naturally
forms a triangle with the pure powers at the vertices, and the conditions of the Theorem are easy to visualize. For example, taking $d=7$, the red sets $N$ in the following pictures
all satisfy the conditions, so the monomials of degree 7 corresponding to the black dots generate ideals with linear presentation:
\begin{center}
\begin{tikzpicture}[scale=.3]
\draw (3.5, 7) circle(.1) [fill=black];
\draw (4, 6) circle(.1) [fill=black];
\draw (3, 6) circle(.1) [fill=black];
\draw (4.5, 5) circle(.1) [fill=black];
\draw (3.5, 5) circle(.1) [fill=black];
\draw (2.5, 5) circle(.1) [fill=black];
\draw (5, 4) circle(.1) [fill=black];
\draw (4, 4) circle(.1) [fill=black];
\draw (3, 4) circle(.1) [fill=black];
\draw (2, 4) circle(.1) [fill=black];
\draw (5.5, 3) circle(.1) [fill=black];
\draw (4.5, 3) circle(.1) [fill=black];
\draw (3.5, 3) circle(.1) [fill=black];
\draw (2.5, 3) circle(.1) [fill=black];
\draw (1.5, 3) circle(.1) [fill=black];
\draw (6, 2) circle(.1) [fill=black];
\draw (5, 2) circle(.1) [fill=black];
\draw (4, 2) circle(.1) [fill=black];
\draw (3, 2) circle(.1) [fill=black];
\draw (2, 2) circle(.1) [fill=black];
\draw (1, 2) circle(.1) [fill=black];
\draw (6.5, 1) circle(.1) [fill=black];
\draw (5.5, 1) circle(.1) [fill=black];
\draw (4.5, 1) circle(.1) [fill=black];
\draw (3.5, 1) circle(.1) [fill=black];
\draw (2.5, 1) circle(.1) [fill=black];
\draw (1.5, 1) circle(.1) [fill=black];
\draw (.5, 1) circle(.1) [fill=black];
\draw (7, 0) circle(.1) [fill=black];
\draw (6, 0) circle(.1) [fill=black];
\draw (5, 0) circle(.1) [fill=black];
\draw (4, 0) circle(.1) [fill=black];
\draw (3, 0) circle(.1) [fill=black];
\draw (2, 0) circle(.1) [fill=black];
\draw (1, 0) circle(.1) [fill=black];
\draw (0, 0) circle(.1) [fill=black];
\draw (3.5, 5) circle(.2) [fill = red];
\draw (4.5, 3) circle(.2) [fill = red];
\draw (5.5, 1) circle(.2) [fill = red];
\draw (3.5, 3) circle(.2) [fill = red];
\draw (4, 2) circle(.2) [fill = red];
\draw (2.5, 3) circle(.2) [fill = red];
\draw (3, 2) circle(.2) [fill = red];
\draw (3.5, 1) circle(.2) [fill = red];
\draw (1.5, 1) circle(.2) [fill = red];
\end{tikzpicture}\vspace{.5cm}

\begin{tikzpicture}[scale=.3]
\draw (3.5, 7) circle(.1) [fill=black];
\draw (4, 6) circle(.1) [fill=black];
\draw (3, 6) circle(.1) [fill=black];
\draw (4.5, 5) circle(.1) [fill=black];
\draw (3.5, 5) circle(.1) [fill=black];
\draw (2.5, 5) circle(.1) [fill=black];
\draw (5, 4) circle(.1) [fill=black];
\draw (4, 4) circle(.1) [fill=black];
\draw (3, 4) circle(.1) [fill=black];
\draw (2, 4) circle(.1) [fill=black];
\draw (5.5, 3) circle(.1) [fill=black];
\draw (4.5, 3) circle(.1) [fill=black];
\draw (3.5, 3) circle(.1) [fill=black];
\draw (2.5, 3) circle(.1) [fill=black];
\draw (1.5, 3) circle(.1) [fill=black];
\draw (6, 2) circle(.1) [fill=black];
\draw (5, 2) circle(.1) [fill=black];
\draw (4, 2) circle(.1) [fill=black];
\draw (3, 2) circle(.1) [fill=black];
\draw (2, 2) circle(.1) [fill=black];
\draw (1, 2) circle(.1) [fill=black];
\draw (6.5, 1) circle(.1) [fill=black];
\draw (5.5, 1) circle(.1) [fill=black];
\draw (4.5, 1) circle(.1) [fill=black];
\draw (3.5, 1) circle(.1) [fill=black];
\draw (2.5, 1) circle(.1) [fill=black];
\draw (1.5, 1) circle(.1) [fill=black];
\draw (.5, 1) circle(.1) [fill=black];
\draw (7, 0) circle(.1) [fill=black];
\draw (6, 0) circle(.1) [fill=black];
\draw (5, 0) circle(.1) [fill=black];
\draw (4, 0) circle(.1) [fill=black];
\draw (3, 0) circle(.1) [fill=black];
\draw (2, 0) circle(.1) [fill=black];
\draw (1, 0) circle(.1) [fill=black];
\draw (0, 0) circle(.1) [fill=black];
\draw (4.5, 3) circle(.2) [fill = red];
\draw (3, 4) circle(.2) [fill = red];
\draw (4.5, 1) circle(.2) [fill = red];
\draw (3, 2) circle(.2) [fill = red];
\draw (2, 2) circle(.2) [fill = red];
\draw (2.5, 1) circle(.2) [fill = red];
\end{tikzpicture}\vspace{.5cm}
\end{center}

In terms of such pictures, the fact that $N$ is a shadow, plus the saturation condition,  means that $N$ is the union of solid upside-down triangles of size 1 or more that do not touch the boundary, while the separation
condition means that the upside-down triangles do not touch one another.
\end{example}

\begin{proof}[Proof of Theorem~\ref{n-1}] Set $S = k[x_1,\dots, x_n]$, and $\mm = (x_1,\dots, x_n)\subset S$.
Given an $\mm$-primary ideal $I$ generated by monomials of degree $d$, we consider the kernel $Y$ of the surjection $S/I \to S/\mm^d$. Set
$\omega = \Ext_S^n(S/\mm^d, S(-n))= \Hom_k(S/\mm^d, k)$ and $\omega_I = \Ext_S^n(S/I, S(-n))= \Hom_k(S/I, k)$. Dual to
$$
0\rTo Y \rTo S/I \rTo S/\mm^d\rTo 0
$$
there is a short exact sequence
$$
0\rTo \omega \rTo \omega_I \rTo X\rTo 0
$$
with $X =  \Ext_S^n(Y, S(-n))= \Hom_k(Y, k)$. We may thus form a (non-minimal) free resolution $\mathbb H$ of $\omega_I$ by the "horse-shoe" construction:
letting $\FF$ and $\GG$ be the minimal free resolutions of $X$ and $\omega$ respectively, the resolution $\mathbb H$ has the form
\begin{small}
$$
F_0\oplus G_0 \lTo^{
\begin{pmatrix}
 d_F&0\\ \phi_1& d_G
\end{pmatrix}
} F_1\oplus G_1 
\lTo \cdots
\lTo F_{i-1}\oplus G_{i-1} \lTo^{
\begin{pmatrix}
 d_F&0\\ \phi_i& d_G
\end{pmatrix}
} F_i\oplus G_i \lTo \cdots
$$
\end{small}
where $d_F$ and $d_G$ are the differentials of $\FF$ and $\GG$, and $\phi_i$ is a map defined inductively: because $F_0$ is
free there is a map
$\epsilon: F_0\to \omega_I$ lifting the augmentation map $F_0 \to X$ along the surjection $\omega_I \to X$; we take $\phi_1$ to be a map lifting $\epsilon$ along the composite
$$
G_0 \to \omega \to \omega_I.
$$
For $i>1$ we let $\phi_i$ be the map lifting the composite $\phi_{i-1}d_{F,i}$ along $d_{G,i-1}$. This construction is summed up in the commutativity of the diagram:
\begin{small}
\begin{center}
\begin{tikzpicture}[scale=.1]
\node(00){$\omega$};
\node(10)[node distance=1cm, below of=00] {$\omega_I$};
\node(20)[node distance=1cm, below of=10] {$X$};

\node(01)[node distance=2cm, right of=00] {$G_0$};
\node(02)[node distance=2cm, right of=01] {$G_1$};
\node(03)[node distance=2cm, right of=02] {$G_2$};
 \node(21)[node distance=2cm, right of=20] {$F_0$};
 \node(22)[node distance=2cm, right of=21] {$F_1$};
 \node(23)[node distance=2cm, right of=22] {$F_2$};
 
 \node(*)[node distance=2cm, left of =10]{$(*)$};

 \node(04)[node distance=7cm, right of=00] {$\cdots$};
\node(14)[node distance=7cm, right of=10] {$\cdots$};
\node(24)[node distance=7cm, right of=20] {$\cdots$};

\draw[->](00) to node{} (10);
\draw[->](10) to node{} (20);
\draw[<-](00) to node{} (01);
\draw[<-](01) to node[above]{$d_{G,1}$} (02);
\draw[<-](02) to node[above]{$d_{G,2}$} (03);
\draw[<-](20) to node{} (21);
\draw[<-](21) to node[below]{$d_{F,1}$} (22);
\draw[<-](22) to node[below]{$d_{F,2}$} (23);

\draw[->,dashed](21) to node[above=3pt, right=0pt]{$\epsilon$} (10);
\draw[->,dashed](22) to node[above=5pt,right= 2pt]{$\phi_1$} (01);
\draw[->,dashed](23) to node[above=5pt,right= 2pt]{$\phi_2$} (02);
\end{tikzpicture}
\end{center}
\end{small}

\begin{lemma}\label{split}
With notation as above, $I$ satisfies $N_{d,p}$ if and only if $\phi_{n-p+1}$ is a split monomorphism, in which case
$\phi_q$ is a split monomorphism for all $q\geq n-p+1$.
\end{lemma}

\begin{proof}[Proof of Lemma~\ref{split}]
 Because $\mm^d$ satisfies $N_{d,n}$, $\omega$ is generated
in degree $-d+1$ and the resolution of $\omega$ is linear except for the last step. Thus $G_i$ is generated in degree $-d+1+ i$ for $i= 0\dots, n-1$. 

On the other hand, 
since $Y = \mm^d/I$ is generated in degree $d$, 
the socle of $X$ is generated in degree $-d$, and thus
$F_n = S(d-n)^{|N|}$, where $|N|$ is the number of monomials of degree $d$ that are not in $I$.
Since
the regularity of $X$ is $-d$,  the generators of $F_i$ have degrees $\leq -d + i$; that is, $F_i$ has the form
$$
 F_i = \bigoplus_j S(d-i+e_{i,j})
$$
with $e_{i,j} \geq 0$.

Because $\FF$ is minimal and the dual of $\FF$ is also acyclic, each $e_ {i,j}$ must be less than or equal to some $e_{i+1,j}$. Thus if,
for a given $i_0$, all the $e_{i_0,j}$ are 0, then $e_{i,j} = 0$ for all $i\geq i_0$ and all $j$.

Now suppose that $\phi_{i_0}$ is a split monomorphism so that, in particular, $e_{i,j} = 0$ for $i\geq i_0$. The 
map $\phi_{i_0}$ takes $\ker d_{F,i_0}$ monomorphically into $\im d_{G,i_0}$. Since $\FF$ is a minimal resolution,
so $\phi_{i_0+1}$ must be a monomorphism. Because all the $e_{i_0+1,j}$ are zero, the free module $F_{i_0+1}$ is
 generated in the same degree as $G_i$, so $\phi_{i_0+1}$ is also a split monomorphism, and repeating the 
 argument we see that $\phi_{i}$ is a split monomorphism for all $i\geq i_0$. Thus the condition that
 $\phi_{n-p_0+1}$ is a split monomorphism is equivalent to the condition that
 $\phi_{n-p+1}$ is a split monomorphism for all $p\leq p_0$.
 
The minimal free resolution of $S/I$ is obtained from $\Hom_S(\mathbb H, S(n))$ by minimizing.
Since $H_i = G_i\oplus F_i$ and the generators of $F_i$ have degree strictly greater than those of $G_i$,
we see that $I$ satisfies $N_{d,p}$ for some $p<n$ if and only if
$e_{i,j} = 0$ for all $i\geq n-p+1$ and $\phi_i$ is a split monomorphism for all $i\geq n-p+1$ so that each summand
$F_i$ in the resolution $\mathbb H$ cancels with a direct summand of $G_{i-1}$ for $i\geq n-p+1$. By the argument
above, this is equivalent  to the condition that $\phi_{n-p+1}$ is itself a split monomorphism.
\end{proof}

To complete the proof of the Theorem we must show that $\phi_1$ is a split monomorphism if and only if
$N$ is the $d$-shadow of a set of monomials whose elements are $(d-1)$-saturated and pairwise $(d-1)$-separated.
First, suppose that $I \subset \mm^d$ satisfies $N_{d,n-1}$ and let $m_1,\dots, m_u$ be a minimal set of monomials
generating the socle of $S/I$, so that $X$ is generated by the dual monomials $\hat m_1,\dots, \hat m_k$.

The socle of $S/\mm^d$ is generated by one monomial of each multi-degree with total degree $d-1$. Thus 
 $\omega$ is generated by one dual monomial of each possible non-negative multi-degree having total degree $-d+1$.
 It follows that for $\phi_1$ to be a split monomorphism, it is necessary that the relations of $X$ contain at most generator
of each multi-degree with total degree $-d+1$. But every monomial of $n$ of degree 
$\deg m_i - d+1$ annihilates $\hat m_i\in X$, so $X$ has relations of multidegree $\deg m_i -\deg n$, and this
must be non-negative; that is, $n$ must divide $m_i$. Thus $m_i$ is $(d-1)$-saturated.

Similarly, if $m_i$ and $m_j$ were not $(d-1)$-separated, then a multiple of each would be equal to the same
monomial of degree $d-1$, and this would give two relations on $X$ with the same multi-degree. This proves that
$m_1,\dots, m_k$ satisfy the conditions of the Theorem, and we must show that $N$ is their shadow.

From the separation condition it follows that the module $X$ is the direct sum of the cyclic submodules $X_i = S\hat m_i$ generated by
the $\hat m_i$. By Lemma~\ref{split}, the module $F_i$ is generated in degree $-d+i$ for all $i\geq 1$, so each $X_i$ has linear resolution
from the first step. It follows that $X_i \cong S/\mm^{\deg m_i - d+1} (\deg m_i)$, and the socle of $X_i$ consists of the duals
of all monomials of degree $d$ that divide $m_i$; thus $N$, which is the union of the duals of
the socles of the $X_i$, is the $d$-shadow.

Conversely, suppose that $N$ is the $d$-shadow of a $(d-1)$-separated set of $(d-1)$ saturated monomials $m_1,\dots, m_u$. It follows
that $m_i\notin I$ but---since $m_i$ is $(d-1)$-saturated, any variable times $m_i$ has a divisor of degree $d$ in $I$, so the $m_i$ generate the part of the socle of $S/I$
of total degree $\geq d$, and the 
$\hat m_i$ generate $X$.
Because the $m_i$ are $(d-1)$-separated, the submodules $X_i\subset X$ intersect in 0, so
$X = \oplus_{1\leq i\leq u} X_i$,  and the socle of $X_i$ is the dual of the $d$-shadow of $m_i$. If $M_i$ is the generator of $F_i$ corresponding
to $m_i$, then 
the first syzygy of $X_i$ is generated by elements $nM_i$  where $n$ ranges over all monomial  of total degree $\deg m_i -d+1$.
Because $m_i$ is $(d-1)$ saturated, $m_i$ is divisible by $n$, and we
 claim that $\phi_1$ may be taken to send $nM_i$ to the generator $N_i$ of $G_0$ that maps to $\widehat{m_i/n}\in \omega$. 
 Since $M_i\in G_0$ maps to $\hat m_i\in X$, the map $\epsilon$ in the diagram $(*)$ may be taken to send $nM_i$ to 
 $\widehat {m_i/n}\in \omega_I$,
 which is also the image of $\widehat{m_i/n} \in \omega$, as required.
\end{proof}
 
 The following result, together with the observation that $Y$ must have a linear resolution
 up to the last step, gives an alternative proof that $X$ is the direct sum of cyclic modules
 of the form $S/\mm^{d_{i}}$ for various $d_{i}$.
 
\begin{proposition}
If $M$ is an indecomposable graded $S$-module of finite length whose first syzygy has linear resolution, then (up to a shift in grading) $M\cong S/\mm^{d}$ for some integer $d$.
\end{proposition}

Graded local duality implies that the socle of $M$ all lies in a single degree and in the cyclic case the result follows---this is the usual (well-known)  proof. In the cyclic case the result also follows from the Herzog-K\"uhl theorem on pure resolutions~\cite{HK}, as in the argument below.

\begin{proof}
Suppose that the generators of $M$ have degrees $g_{1}, \dots, g_{t}$, and that the relations are all in degree $d$. We will show that $M \cong \oplus_{i} S(-g_{i})/\mm^{d-g_{i}}$.

 Let $P$ be the minimal presentation matrix of $M$, with
$i$-th row corresponding to a generator of degree $g_{i}$. Set $s_{i} = {n-1+(d-g_{i})\choose n-1}$. By Boij-S\"oderberg theory \cite[Theorem 0.2]{ES}, the Betti table of $M$ is the sum of the Betti tables of
the modules $S(-g_{i})/\mm^{d-g_{i}}$, and in particular $P$ has $\sum_{i}s_{i}$ columns. If the forms of degree $d-g_{i}$ in the i-th column of $P$ span a space of dimension $e_{i} \leq s_{i}$, then
after suitable column transformations $P$ would have $\sum_{i}s_{i} - \sum_{i}e_{i}$ columns
of zeros, so $e_{i} = s_{i}$ for all $i$. In the case $t=1$ the result follows immediately. (Note that this case does not require the full force of \cite[Theorem 0.2]{ES}, since when $t=1$ the resolution is pure, and its shape is given by the Herzog-K\"uhl Theorem).

We will prove by induction on $t$ that any $t\times \sum_{i}s_{i}$ matrix $P$ without columns of zeros, whose $i$-th row contains forms of degree $d-g_{i}$, and whose maximal minors generate an
$\mm$-primary ideal, is the direct sum of 1-rowed matrices such as 
$$\begin{pmatrix}
f_{1}&\cdots &f_{s_{1}}&0&\dots&0&0&\dots&\dots&0\\
0&\dots&0&f'_{1}&\cdots &f'_{s_{2}}&0&\dots&\dots&0\\
0&\dots&0&0&\dots&0&f_{1}''&\dots
\end{pmatrix}.
$$

The case $t=1$ is trivial. After a suitable column operation we may  assume that the first $s_{1}$ columns of the first
row of $P$ contain a basis for the forms of degree $d-g_{1}$, and that all the other entries of the first row are zero. Let $P'$ be the submatrix of $P$ omitting the first row and 
the first $s_{1}$ columns. The ideal of maximal minors of $P$ is the $\mm^{d-g_{1}}$ times the ideal of maximal minors
of $P'$, so $P'$ satisfies the same hypotheses as $P$. By induction $P'$ has the desired form. After further column operations we may assume that the
first $s_{1}$ columns of $P$ have zeros in all but the first row, and thus
$P$ is equivalent to a direct sum of 1-rowed matrices. From this it follows that
$M$ is a direct sum of cyclic modules, so $M$ must be cyclic, and we are done.
\end{proof}

\section{Examples: fractal ideals and obstructions to shelling}\label{sec_last}

In this section we use our previous results to construct examples of $N_{d,p}$ $\mm$-primary ideals with interesting behavior. 
Our first construction is inspired by fractal geometry. It gives us ideals with relatively few generators that have almost linear resolution. We start with the case $n=3$, where a pattern is easiest to describe. 

\begin{proposition}\label{sier} With $n=3$, 
for $d = 1,2$ set:
$I_d=\mm^d$. Inductively, define:\\
$I_{2r-1} = (x_1^r,x_2^r,x_3^r)I_{r-1}$,\\
$I_{2r} = x_1^{r+1}I_{r-1} + (x_2^r, x_3^r)I_r$.\\
The ideal $I_d$ has linear presentation for all $d$. 
\end{proposition}

\begin{proof}
In terms of diagrams as in Example \ref{ex3}, the ideals $I_d$ are created by staring from the simplex of monomials of degree $d$,  and first removing the largest upside-down triange of monomials that does not meet the boundary. This leaves three smaller simplices, and we repeat the pattern within each of them, etc. The result for $d=7$ is
shown in the first diagram of Example \ref{ex3}. From Theorem \ref{n-1} we see that this produces an ideal satisfying $N_{d,n-1}$.
\end{proof}

\begin{remark}\label{rem3}
When $d=2^r-1$,  $I_d= \mm\mm^{[2]}\mm^{[4]}\dots\mm^{[2^{r-1}]}$. As the reader may show, this has precisely $3^r$ minimal generators which is $O(d^{log_23})$. By contrast, $\mm^d$ has $O(d^2)$ generators. 

In this situation, the picture of the generators of $I_d$ is exactly the so-called Sierpi\'nski triangle or gasket. See Figure 1.  We note that (generalizations of) Sierpi\'nski triangles have also appeared in \cite{FMMS} where they were used to compute Frobenius powers of monomial ideals. 
\begin{figure}
  \centering
  \includegraphics{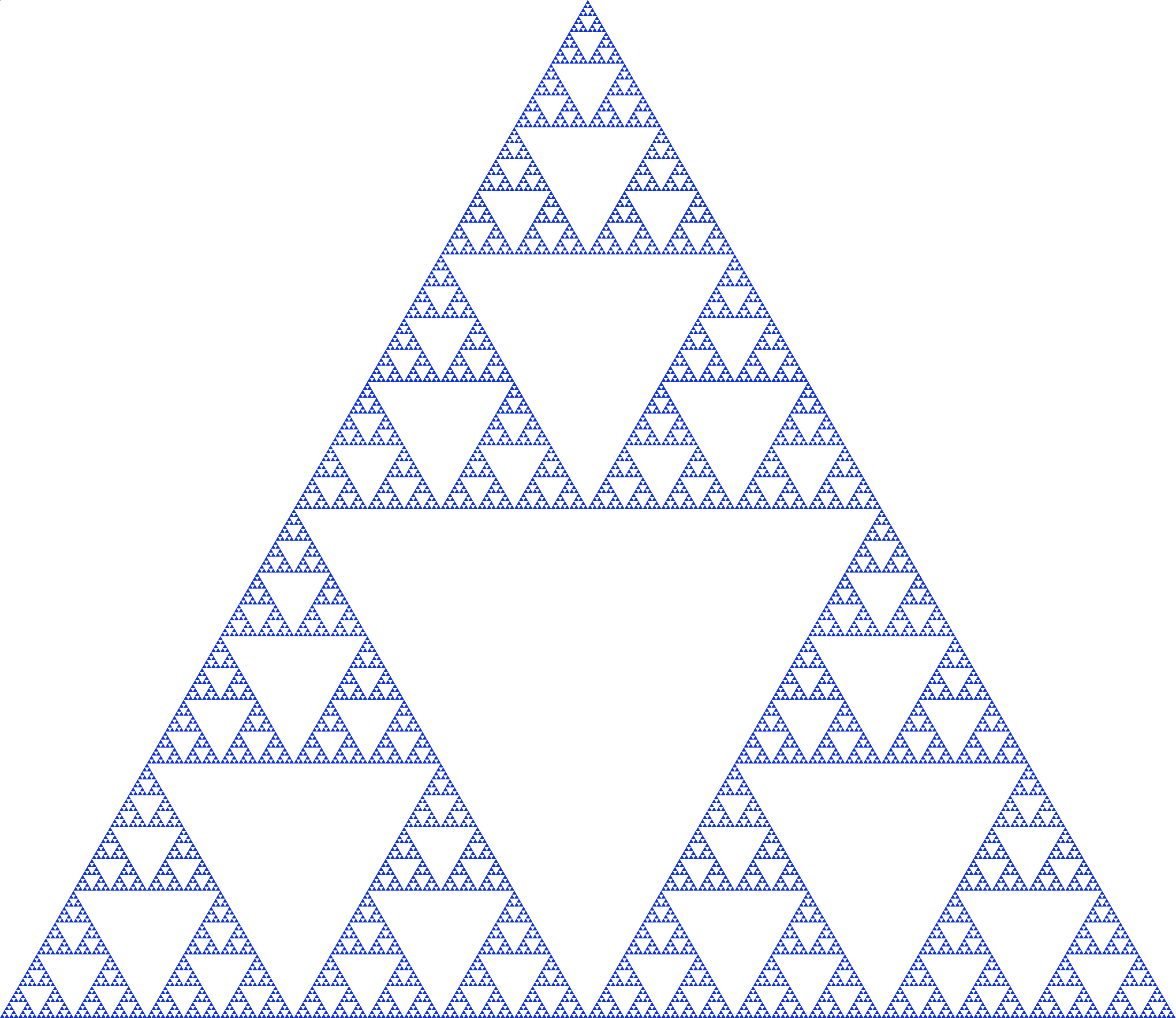}
  \caption{Sierpi\'nsky triangle by Beojan Stanislaus, CC BY-SA 3.0, https://commons.wikimedia.org/w/index.php?curid=8862246}
\end{figure}

\end{remark}

For  $n>3$ variables, it is harder to give a recursive formula for such sparse ideals with $N_{d,p}$. Instead we offer a closed form formula for special values of $d$:

\begin{proposition}\label{sier>3}
For any $n$, each of the following ideals
$$
I_r =  (\mm\mm^{[2]}\mm^{[4]}\dots\mm^{[2^{r-1}]})^{p-1}\mm^{[2^r]} \subset S
$$
satisfies $N_{d,p}$.  
\end{proposition}

\begin{proof}
We do induction on $r$, and first consider $I_1 = \mm^{p-1}\mm^{[2]}$. Since $t_{p-1}(\mm^{[2]})= 2p$, $\mm^{[2]}\mm^{p-1} = \mm^{[2]}\cap \mm^{p+1}$ has  linear resolution for $p-1$ steps by Proposition~\ref{trunc}.  

For the induction step, write $I_{r+1}$ as $\mm^{p-1}J$, and note that $J$ is the ideal $I_r$ constructed inside  the polynomial ring $S'=k[x_1^2,\dots, x_n^2]$. By induction, $J$ has a resolution that is $p-1$ step linear in $S'$, which means that over $S$, $J$ has a resolution which is quadratic in $p-1$ steps. That implies $t{p-1}(J)= 2p$, and hence $I_{r+1} =\mm^{p-1}J$ satisfies $N_{d,p}$, as desired. 
\end{proof}

\begin{remark}\label{rem>3}
As in \ref{rem3}, it can be shown that the ideal $I_{r}$ constructed in \ref{sier>3} are much sparser than $\mm^d$, even in the almost linear resolution case, $p=n-1$.
\end{remark}

Next, we discuss an application to shelling. We begin by  defining the algebraic analogue of shelling for monomial ideals. 

\begin{definition}
Let $I$ be a monomial ideal generated in degree $d$ and $f$ is a degree $d$ monomial. We say that the transition $I\mapsto (I,f)$ is a shelling move if $I:f$ is generated by a subset of variables. We say that an ideal $L$ is shelled over $I$ if it can be obtained from $I$ by a sequence of shelling moves. 
\end{definition}

The next result is the algebraic version of \cite[Lemma 3.1]{DDL}.

\begin{proposition}\label{shellprop}
Let $I$ be a monomial ideal generated in degree $d$ and $f$ be a monomial of degree $d$. 
\begin{enumerate}
\item If $(I,f)$ is $N_{d,2}$ then $I\mapsto (I,f)$ is a shelling move. 
\item If $I$ satisfies $N_{d,2}$, then  $I\mapsto (I,f)$ is a shelling move if and only if $(I,f)$ also satisfies $N_{d,2}$.
\item If $I$ satisfies $N_{d,p}$ for some $p\geq 2$ and $I\mapsto (I,f)$ is a shelling move, then $(I,f)$ also satisfies $N_{d,p}$. 
\end{enumerate}
\end{proposition}

\begin{proof}
We have a short exact sequence of graded $S$-modules
\[0 \to S/(I:f)(-d) \to S/I \to S/(I,f)) \to 0.\]
This exact sequence induces the following long exact sequence in $\Tor$: 
\begin{align*}
 \cdots \to \Tor_2(S/(I,f)),k) \to &\Tor_1^S(S/(I:f)(-d),k) \to\\
   &\Tor^S_1(S/I,k) \to \Tor^S_1(S/(I,f)),k)\to \cdots
\end{align*}

The map $\Tor^S_1(S/I,k) \to \Tor^S_1(S/(I,f)),k)$ is injective. Hence the map $\Tor_2(S/(I,f),k) \to \Tor_1^S(S/(I:f)(-d),k)$ is surjective.  But since $S/(I,f)$ has linear first syzygy so $\Tor_2^S(S/(I,f),k) \cong k^{\beta^S_1(I,f)}(-d-1)$. Hence $\Tor_1^S(S/(I:f)(-d),k)$ is generated in degree $-d-1$.  Thus $(I:f)$ is generated in degree $1$. The proof of $(b)$ and $(c)$ are similar.

\end{proof}


We can characterize when an $N_{d,n-1}$ ideal is shelled over another. For an monomial ideal generated in degree $d$ we write $N(I)$ for the set of monomials of degree $d$ not in $I$

\begin{corollary}\label{noshell}
If $I,J$ are $\mm$-primary monomial ideals satisfying $N_{d,n-1}$, then $J$ is shelled over $I$ if and only if $N(I)$ is a disjoint union of $N(J)$ and a set of singleton shadows $M=\{m_1,\dots, m_s\}$. If that is the case, the shelling can be obtained by adding elements in $M$ in any order. 
\end{corollary}

\begin{proof}
If there were a shelling from $I$ to $J$ then each intermediate ideal would also satisfy $N_{d,n-1}$ by  \ref{shellprop}. But by our structure Theorem \ref{n-1}, the difference between $N(I)$ and $N(J)$ is a disjoint unions of $d-1$ saturated shadows. But the $d-1$-saturated shadow of a monomial is a simplex in the monomial lattice, so if it is not a singleton, then after removing one monomial it is no longer a shadow.  

On the other hand, if $N(J)\setminus N(I)$ is a union of singleton shadows, then one can fill them in one by one to get from $I$ to $J$ in any order, and such collections of moves are shelling as each intermediate ideal is $N_{d,n-1}$ by  Theorem \ref{n-1} and  \ref{shellprop}(2). 
\end{proof}

This result implies a rigidity of regularity:

\begin{corollary}\label{rigidity}
Supose that $I$ is a primary monomial ideal satisying $N_{d,n-1}$. If $\reg I\geq d+2$ and $I\mapsto (I,f)$ is a shelling move, then $\reg (I,f) = \reg I$.
\end{corollary}

\begin{proof}
The larger simplices in $N(I)$ cannot be changed by a shelling move.
\end{proof}

\begin{example}\label{ex_noshell}
Consider the ideal $I=I_6= (x^4)\mm^2+(y^3,z^3)\mm\mm^{[2]}$ in \ref{sier}.  The set $N(I)$ contains the shadow of $x^3y^2z^2$, which is the triangle $\{x^3y^2z,x^3yz^2,x^2y^2z^2\}$. It follows that $\mm^6$ is not shelled over $I$. 
\end{example}

From Corollary \ref{noshell} we can deduce a similar result for square-free monomial ideals or, equivalently, simplicial complexes. If $I$ is generated by monomials
of degree $d$, we let $T$ be the polynomial ring $k[x_{ij}]_{1\leq i\leq n, 1\leq j\leq d}$. Define $\pol: S \to T$ to be the map on monomials that takes $x_i^r$ to $x_{i1}x_{i2}\dots x_{ir}$, and let $\depol: T \to S$ be the map of algebras that takes $x_{ij}$ to $x_i$. Note that polarization commutes with $\lcm$, and that
$\depol (\pol(m)) = m$.

\begin{lemma}\label{depol}
If $J\subset T$ is a square-free monomial idea satisfying $N_{d,2}$ then $I := \depol(J)$ satisfies $N_{d,2}$ as well.
\end{lemma}

\begin{proof}
By Proposition~\ref{lppath}, it suffices to show that given any two minimal generators $f = \depol(F), g = \depol(G)$ of $I$ are connected by a path within the support of $\lcm(f,g)$. Since $J$ satisfies $N_{d,2}$,  Proposition~\ref{lppath} shows that there is a path  of monomials in $J$ of monomials dividing $\lcm(F,G)$. Depolarizing these monomials we get a path in support of $\lcm(f,g)$.
\end{proof}


\begin{corollary}\label{cor_depol}
Let $I,J$ be $N_{d,p}$ square-free monomial ideals in $T$ for some $p\geq 2$. If $J$ is shelled over $I$, then $\depol(J)$ is shelled over $\depol(I)$. 

Consequently, if $I,J$ are monomial ideals in $S$ and $J$ is not shelled over $I$, then $\pol(J)$ is not shelled over $\pol(I)$. 
\end{corollary}

\begin{proof}
Let $J_0=I\mapsto \cdots \mapsto  J_s=J$ be a sequence of shelling moves from $I$ to $J$. By Lemma \ref{depol}, each ideal $\depol(J_I)$ is $N_{d,2}$. Proposition \ref{shellprop} implies that $\depol(J_1)\mapsto \cdots \mapsto \depol(J_s)$ is a sequence of selling moves, so $\depol(J)$ is shelled over $\depol(I)$. 
The second statement follows from the first and the fact that $\depol(\pol(I))=I$. 
\end{proof}

\begin{example}\label{simp_noshell}
Let $T$ be the ring $k[x_{ij}]_{1\leq i\leq 3, 1\leq j\leq 6}$ and $S=k[x_1,x_2,x_3]$. The ideal $\depol(\m^6)$ is not shelled over $J=\depol(I)$, where $I$ is the ideal in Example \ref{ex_noshell}. 

Taking the Alexander dual, one obtains a normal (equivalently, satisfying Serre's condition $(S_2)$) simplicial complex $\Delta = \Delta(J^\vee)$ of dimension $11$ such that the complete $11$-skeleton of the $17$-simplex is not shelled over $\Delta$. 

Of course, if we apply the same process using \ref{noshell} with $n\geq 3$ variables one gets simplicial complexes $\Delta$ satisfying Serre's condition $S_{n-1}$ that can not be extended by shelling to the complete $\dim \Delta$-skeleton of the simplex on all vertices. That is because the polarization of $N_{d,p}$ ideals are also $N_{d,p}$, and the Alexander dual of $N_{d,p}$ ideals define simplicial complexes that satisfy Serre's condition $(S_p)$, see \cite{MT}. 
\end{example}

\section{Questions and discussions}\label{sec_ques}

In this final section we collect some questions inspired by the literature and our own work. Let $\mathcal C$ be a class of monomial ideals in $S$. The most prominent examples are $\mathcal C=\{\text{square-free ideals}\}$ or $\mathcal C=\{\text{primary ideals}\}$. Let $N_{d,p}(\mathcal C)$ denote the ideals in $\mathcal C$ that are $N_{d,p}$. 

\begin{question}\label{q1} 
What can we say about the Betti numbers of ideals in $N_{d,p}(\mathcal C)$? For instance, it is intuitively clear that such ideals must not have too few generators. Can we prove good bounds? What about optimal examples? 
\end{question}

The only result we are aware of in this direction is \cite[Proposition 11.1]{EHU}, where it was proved that any graded $\m$-primary ideal with almost linear resolution must have at least $\binom{n+d-2}{d}+\binom{n+d-3}{d-1}$ generators, with equality if and only if $S/I$ is Gorenstein. See \cite{BG} for a recent survey of the literature on lower bounds for Betti numbers  of ideals in general and monomial ideals in particular.

Equally sensible is the expectation that ideals in $N_{d,p}(\mathcal C)$ must have low Castelnuovo-Mumford regularity. 

\begin{question}\label{q2} 
Can we establish sharp upper bounds for $\reg I$, $I\in  N_{d,p}(\mathcal C)$? What about optimal examples? 
\end{question}

In this direction, there is a $O(\log(n))$ bound on regularity of  $N_{2,2}$ monomial ideals (using \cite{DHS} for the square-free case and polarization). Interestingly, we only know monomial $N_{2,2}$ ideals with regularity $O(\log(\log(n))$, using constructions from the study of hyperbolic Coxeter groups (\cite{CKV}). For $d>2$, it has been conjectured that square-free $N_{d,2}$ ideals have regularity at most  $n - \lfloor \frac{n}{d+1} \rfloor - \lfloor \frac{n-1}{d+1} \rfloor$ (\cite{DT, DV}). Only the case $d=3$ has been settled (\cite{DV}). 

One can sometimes show that an ideal is in $\mathcal C$ is $N_{d,p}$ by checking the restrictions to all subsets of variables of size $r$, for some relatively small value of $r$,
as in Corollary~\ref{local} and Theorem \ref{primecubic}. If that is the case we say that $N_{d,p}(\mathcal C)$ is $r$-certifiable. For instance, if $\mathcal C$ is the class of all monomial ideals, $N_{d,p}(\mathcal C)$ is $dp$-certifiable. In the quadratic case, much better bound is known, indeed, $N_{2,p}(\mathcal C)$ is $(p+2)$-certifiable. Note that these bounds do {\bf not}  depend on $n$.

\begin{question}\label{q3} 
Given $d,p$, what is the smallest value of $r$ such that $N_{d,p}(\mathcal C)$ is $r$-certifiable? 
\end{question}

\bibliographystyle{ABC99}

\begin{thebibliography}{ABC99}


\bibitem{BHZ}[BHZ], M. Bigdeli, J. Herzog and R. Zaare-Nahandi, 
\emph{On the index of powers of edge ideals},
Commun. Algebra 46 (2018), No. 3, 1080--1095.
 

\bibitem {B}[B] A. Boocher,
\emph{Free resolutions and sparse determinantal ideals}, 
Math. Res. Lett. 19 (2012) 805–821. 

\bibitem{BG}[BG] A. Boocher and E. Grifo, \emph{Lower bounds on Betti numbers}, preprint,  \url{https://arxiv.org/abs/2108.05871}, 2021. 

\bibitem{CKV}[CKV] A. Constantinescu, T. Kahle, and M. Varbaro, \emph{Linear syzygies, flag
complexes, and regularity}, Collect. Math. 67 (2016), no. 3, 357--362. 

\bibitem{DHS}[DHS] H. Dao, C. Huneke and J. Schweig, \emph{Bounds on the regularity and projective
dimension of ideals associated to graphs}, J. Algebraic Combin. 38 (2013), no. 1, 37--55.

\bibitem{DT}[DT] H. Dao and S. Takagi, \emph{On the relationship between depth and cohomological dimension}, Compositio Math. 152 (2016), no.4, 876--888.

\bibitem{DV}[DV]  H. Dao and T. Vu, \emph{Regularity of  monomial ideals with  linear syzygies}, preprint. 

\bibitem{DDL}[DDL] H. Dao, J. Doolittle and J. Lyle, \emph{Minimal Cohen-Macaulay simplicial complexes}, SIAM J. Discrete Math. 34 : 3(2020), 1602--1608.

\bibitem{FMMS}[FMMS] C. Francisco, M. Mastroeni, J. Mermin and J. Schweig, 
\emph{Computing Frobenius powers of monomial ideals}, preprint. 


\bibitem{GPW}[GPW], V. Gasharov, I. Peeva and V. Welker  \emph{The lcm lattice in monomial resolutions}, Math. Res. Lett. 6
(1999), no. 5-6, 521--532.


\bibitem{EGHP}[EGHP] D. Eisenbud, M. Green, K. Hulek and S. Popescu, \emph{Restricting linear syzygies: algebra and geometry},
Compositio Math. 141 (2005) 1460--1478.

\bibitem{EHU}[EHU]  D. Eisenbud, C. Huneke and B. Ulrich, \emph{Regularity of Tor and graded Betti numbers}, Amer. J. Math. 128, 3 (2006), 573--605.


\bibitem{ES}[ES] D. Eisenbud and F.-O. Schreyer, \emph{Betti numbers of graded modules and cohomology of vector bundles}, Journal of the American Mathematical Society 22 (2009), no. 3, 859--888.

\bibitem{FSY}[FSY] M. Farrokhi, Y. Sadegh and A. Yazdan Pour, \emph{Green-Lazarsfeld index of square-free monomial ideals and their powers}, preprint, \url{https://arxiv.org/abs/2110.12174}, 2021.


\bibitem{Fr}[Fr] R. Fr\"oberg, \emph{On Stanley-Reisner rings}, Topics in Algebra, Part 2 (Warsaw, 1988), Banach Center Pub. 26 (1990) PWN Warsaw.

\bibitem{HK}[HK] J. Herzog and M. K\"uhl, \emph{On the Betti numbers of finite pure and linear resolutions}, Comm. Algebra 12 (1984) 1627--1646.

\bibitem{MN}[MN] J. Migliore and U. Nagel, \emph{A Tour Of The Weak And Strong Lefschetz Properties}, J. Comm. Alg., 5
(2013), 329--358.

\bibitem{MT}[MT] S. Murai and T. Terai, \emph{$h$-vectors of simplicial complexes with Serre's $S_c$ conditions}, Math. Res. Lett. {\bf 16} (2009), 1015--1028.


\bibitem{PV}[PV] I. Peeva and M. Velasco,  \emph{Frames and degenerations of monomial resolutions}, Trans. Amer. Math. Soc. 363, no. 4  (2011), 2029--2046.

\end{thebibliography}

\end{document}